\documentclass[12pt,leqno]{article}
\usepackage{amsopn,amsmath,amssymb,amscd,latexsym}
\usepackage[all]{xy}
\newcommand{\C}{{\mathbb{C}}}
\newcommand{\F}{{\mathbb{F}}}
\newcommand{\Q}{{\mathbb{Q}}}
\newcommand{\Z}{{\mathbb{Z}}}
\newtheorem{theorem}{Theorem}%[section]
\newtheorem{lemma}[theorem]{Lemma}
\newtheorem{prop}[theorem]{Proposition}
\newenvironment{proof}{\noindent {\bf Proof}}{\mbox{}\hfill$\Box$}
\begin{document}
%%%%%%%%%%%%%%%%%%%%%%%%%%%%%%%%%%%%%
\title{On the irreducibility of the two variable zeta-function for
curves over finite fields}
%\author{Christopher \surname{Deninger} \email{deninge@math.uni-muenster.de}}
\author{Niko Naumann}
%\institute{Mathematisches Institut, WWU M\"unster, Einsteinstr. 62, 48149 M\"unster, Germany}
%\author{Wilhelm \surname{Singhof} \email{singhof@cs.uni-duesseldorf.de}}
%\institute{Mathematisches Institut, Universit\"at D\"usseldorf, Universit\"atsstr. 1, 40225 D\"usseldorf, Germany}
%\runningauthor{Christopher Deninger, Wilhelm Singhof}
%\runningtitle{Real polarizable Hodge structures arising from foliations}
\date{\today }
\maketitle
\begin{abstract}
  In \cite{P} R. Pellikaan introduced a two variable zeta-function $Z(t,u)$
for a curve over a finite field $\F_q$ which, for $u=q$, specializes to the usual zeta-function and he proved, among other things, rationality:
$Z(t,u)=(1-t)^{-1}(1-ut)^{-1}P(t,u)$ with $P(t,u)\in\Z [t,u]$. We prove that $P(t,u)$
is absolutely irreducible. This is motivated by a question of J. Lagarias and
E. Rains about an analogous two variable zeta-function for number fields.
\end{abstract}
%%%%%%%%%%%%%%%%%%%%%%%%%%%%%%%%%%%%%


\section{Introduction}
\label{sec:1}
Let $X$ be a proper, smooth, geometrically connected curve of genus $g$ over
the finite field $\F_q$. The zeta-function of $X/\F_q$ can be written as a power series
\[
Z(t)=\sum_D\frac{q^{h^0(D)}-1}{q-1}t^{deg(D)}.
\]
Here the sum is over $\F_q$-rational divisor classes of $X$ and\\
$h^0$($D$):= dim$_{\F_q}$ $H^0(X,{\cal O}(D))$. Writing $b_{nk}$ for
the number of divisor classes of degree $n$ and with $h^0(D)$=$k$ this becomes
\[
Z(t)=\sum_{n\ge 0}\sum_{k\ge 1} b_{nk}\frac{q^k-1}{q-1}t^n.
\]
In \cite{P} R. Pellikaan observed that the classical proof of rationality and the functional equation for $Z(t)$ go through when $q$ is treated as a variable in this expression. He thus introduced the following power series in \cite{P}, Def. 3.1:

\[
Z(t,u):=\sum_{n\ge 0}\sum_{k\ge 1} b_{nk}\frac{u^k-1}{u-1}t^n.
\]

This is called the two variable zeta-function of the curve.
We will denote by $h$ the class-number of $X/\F_q$, i.e. $h=$$|$Pic$^0$($X$)$|$. Then Pellikaan proved:

\begin{theorem}\label{rationality}

We have
\[
Z(t,u)=(1-t)^{-1}(1-ut)^{-1}P(t,u) \mbox{ with } P\in\Z [t,u].
\]
Furthermore:\\
1) deg$_t$ $P$=$2g$, deg$_u$ $P$=$g$.\\
2) In the expansion $P(t,u)=\sum_{i=0}^{2g} P_i(u)t^i$ one has
$P_0(u)=1$,\\
deg$_u$ $P_i(u)\leq i/2+1$ and $P_{2g-i}(u)=u^{g-i}P_i(u)$ for
$0\leq i\leq 2g$.\\
3) $P(1,u)=h$.

\end{theorem}

Here deg$_u$ and deg$_t$ denote the degree of a polynomial in the indicated variable.
The above results are all taken from \cite{P}, Prop. 3.5 and we copied only
those needed later on. Note that the statement deg$_u$ $P_i(u)\leq i/2$
in [loc. cit.] is a misprint.
Indeed, we will see below that one always has deg$_u$ $P_1(u)=1$ (unless
$g=0$). As expected we have $P(t,u)=1$ in case $g=0$.\\

In \cite{GS} G. van der Geer and R. Schoof used analogies from
Arakelov-theory to define a two variable zeta-function for number fields
along the above lines. As the number field case will serve only as a motivation
in this note we refer to the original sources \cite{GS} and \cite{LR} for definitions
and to \cite{D} for a comparison between them. Suffice it to say
that in \cite{LR}, section 8 we find an entire function $\xi_{\Q}(w,s)$ of two
complex variables which for $w=1$ equals Riemann's $\xi$-function.
In particular, the zeroes of $\xi_{\Q}(1,s)$ are precisely the non-trivial
zeros of the Riemann zeta-function. One is thus led to study the
zero-locus of $\xi_{\Q}(w,s)$. J. Lagarias and E. Rains ask whether
it might be the closure of a single irreducible complex-analytic
variety of multiplicity one. The corresponding question in the
geometric case seems to be whether the zero-locus of $P(t,u)$ is
irreducible.\\
This is indeed the case:

\begin{theorem}\label{irr}

In the above situation, $P(t,u)$ is irreducible in $\C(u)[t]$.

\end{theorem}

As an illustration we discuss the cases $g=1$ and $g=2$:\\
For $g=1$ setting $N:=|X(\F_q)|$
 we have
\[
P(t,u)=1+(N-1-u)t+ut^2,
\]
c.f. \cite{P}, example 3.4. This polynomial is reducible in $\C(u)[t]$
if and only if $N=0$ in which case we have $P(t,u)=(1-t)(1-ut)$. But
it is well-known that a curve of genus one over a finite field
 always has a rational point, i.e. $N\neq 0$.\\
In case $g=2$ let the usual zeta-function of $X$ be
\[
Z(t)=(1-t)^{-1}(1-qt)^{-1}L(t)
\]
with
\[
L(t)=1+at+bt^2+qat^3+q^2t^4
\]
for certain $a,b\in\Z$. As $X$ is hyperelliptic, Prop. 4.3. of
\cite{P} can be used to compute
\[
P(t,u)=1+((a+q)-u)t+((q(q-1)+aq+b)-(a+q-1)u)t^2+((a+q)-u)ut^3+u^2t^4.
\]
This will not be used in the sequel and we omit the proof.\\
In order not to lead intuition astray we point out that in general
$Z(t,u)$ is not determined by $Z(t)$, see \cite{P}, example 4.4.\\
After a lengthy computation with discriminants one sees that a neccessary
condition for this $P(t,u)$ to be reducible is
\[
b+a(q+1)+(q^2+1)=0.
\]
But this expression equals $L(1)=h\neq 0$ !\\
The fact that $h\neq 0$ enters in the general proof of theorem \ref{irr}
precisely through the condition $\beta\neq 0$ of lemma \ref{irrcrit} below.
Note, however, that condition 2) of this lemma cannot be dropped. So one
needs one more result on $P(t,u)$, contained in Prop. \ref{leadingcoeff},
which follows from Clifford's theorem.\\

I would like to thank C. Deninger for posing the above problem and for useful discussions
on the topic.

\section{Proof}

We will use the following criterion for irreducibility:

\begin{lemma}\label{irrcrit}

Let $k$ be a field, $F\in k[u,t]$ and assume:\\
1) $F$ is monic in $t$.\\
2) the leading coefficient of $F$ as a polynomial in $u$ is
irreducible in $k[t]$.\\
3) there are $\alpha,\beta\in k$, $\beta\neq 0$ with $F(u,\alpha)=\beta$.

Then $F$ is irreducible in $k(u)[t]$.

\end{lemma}

This lemma will be applied to
\[
F=\tilde{P}(t,u):= t^{2g}P(t^{-1},u)\in\C [u,t].
\]
Note that the irreducibility of $\tilde{P}$ in $\C(u)[t]$ will imply
the irreducibility of $P$ because $P(0,u)=P_0(u)=1\neq 0$ by theorem \ref{rationality}.
The advantage of $\tilde{P}$ is that it is monic in $t$ and so satisfies condition
1) of lemma \ref{irrcrit}. Also 3) is satisfied (with $\alpha=1$, $\beta=h$) according to theorem \ref{rationality}, 3).

\begin{proof} (of lemma \ref{irrcrit}) Assume to the contrary that $F=fg$
in $k(u)[t]$ with $f$ and $g$ of positive degree and monic. One knows, c.f. for example \cite{E}, Prop 4.11, that the coefficient of $f$ and $g$ are integral over $k[u]$ and as $k[u]$ is integrally closed we have $f,g\in k[u,t]$.
So we can consider the decomposition $F=fg$ as polynomials in $u$ and infer
from 2) that the leading coefficient of $f$ as a polynomial in $u$ lies in $k[t]^*=k^*$ (upon
exchanging $f$ and $g$ if neccessary). In particular $n:=$deg$_u$ $f(u,t)=$deg$_u$ $f(u,\alpha)$. Substituting 3) gives $\beta=f(u,\alpha)g(u,\alpha)$ in
$k[u]$. As $\beta \neq 0$ we get $n=0$, i.e. $f$ is constant in $u$ hence
$f\in k^*$, contradiction.

\end{proof}

We are left with verifying condition 2) of lemma \ref{irrcrit} for
the given $\tilde{P}$, i.e. the leading coefficient of $\tilde{P}$ as
a polynomial in $u$ is irreducible in $k[t]$. We will in fact determine
this coefficient:

\begin{prop}\label{leadingcoeff}

For $g\geq 1$: $\tilde{P}(t,u)=(1-t)u^g+O(u^{g-1})$.

\end{prop}

\begin{proof} We already know deg$_u$ $\tilde{P}=g$. Also the
assertion is clear for $g=1$ from the formula for $P(t,u)$ recalled
 in the introduction. We assume $g\ge 2$. Looking at
\[
\tilde{P}(t,u)= t^{2g}+P_1(u)t^{2g-1}+\cdots +P_g(u)t^g+uP_{g-1}(u)t^{g-1}+\cdots +u^gt^0
\]
and using the bound deg$_u$ $P_i(u)\leq i/2+1$ we see that $u^g$ can
only occur in the last three terms: $u^{g-2}P_2(u)t^2+u^{g-1}P_1(u)t+u^g$.
So the proof is completed by the following result on $P_1$ and $P_2$.

\end{proof}

\begin{prop} For $g\geq 1$:\\
1) deg$_u$ $P_1(u)=1$ and the leading coefficient is $-1$.\\
2) deg$_u$  $P_2(u)\leq 1$.

\end{prop}

\begin{proof} This is again clear for $g=1$. We assume $g\ge 2$ and
write $P_i(u)=\sum_k\alpha_{ik}u^k$, $\alpha_{ik}\in\Z$.
As we already know deg$_u$ $P_1(u)\leq 1$ and deg$_u$ $P_2(u)\leq 2$
we need to show $\alpha_{11}=-1$ and $\alpha_{22}=0$. Recalling the
notation $b_{nk}$ from the introduction we have the following

{\bf Claim}: $b_{12}=\alpha_{00}+\alpha_{11}$ and $b_{23}=\alpha_{00}+\alpha_{11}+\alpha_{22}$.

Granting this we observe that Clifford's thorem, c.f. \cite{H}, IV, thm 5.4, gives $b_{12}=b_{23}=0$. Recalling also $\alpha_{00}=1$, because $P_0(u)=1$, and substituting gives indeed $\alpha_{11}=-1$ and $\alpha_{22}=0$.

\end{proof}

To prove the above claim we write the rational expression for $Z(t,u)$ in theorem \ref{rationality} in terms
of coefficients:

\[
Z(t,u)=\sum_{n\geq 0}\sum_{k\geq 1} b_{nk}\frac{u^k-1}{u-1}t^n=(\sum_{i\geq 0}t^i)(\sum_{j\geq 0} (ut)^j)(\sum_{l\geq 0}\sum_{k\geq 0} \alpha_{lk} u^kt^l).
\]

This gives

\begin{lemma}\quad\\
1) for $\nu$,$\alpha\geq 0$ : $\sum_{k\geq\alpha +1} b_{\nu k}=\sum_{\mu,i\geq 0,\mu+i\leq\nu}\alpha_{i,\alpha-\mu}$\\
2) for $\nu\geq 0$, $\mu\geq 1$: $b_{\nu\mu}=\sum_{i=0}^{\nu}(\alpha_{i,\mu-\nu-1+i}-\alpha_{i\mu})$

\end{lemma}

We omit the details of this straightforward computation except to say
that for 1) one uses $\frac{u^k-1}{u-1}=1+\cdots+u^{k-1}$ and 2) follows
from 1) by a telescope-summation. In the formulation of the lemma it is
understood that $\alpha_{nk}=0$ whenever $k<0$.
We get from 2):\\
\[
b_{23}=\alpha_{00}-\alpha_{03}+\alpha_{11}-\alpha_{13}+\alpha_{22}-\alpha_{23}
\]
and
\[
b_{12}=\alpha_{00}-\alpha_{02}+\alpha_{11}-\alpha_{12}.
\]
But we know $\alpha_{02}=\alpha_{03}=\alpha_{12}=\alpha_{13}=\alpha_{23}=0$ because deg$_u$ P$_i$(u)$\leq i/2+1$.\\
This concludes the proof of the claim, hence of theorem \ref{irr}.

Mathematisches Institut der WWU M\"unster\\
Einsteinstr. 62\\
48149 M\"unster\\
Germany\\
e-mail: naumannn@uni-muenster.de\\

\vspace*{0.5cm}

\hspace*{\fill} \\
\end{document}